--------------------------------------------------------------------------
\documentstyle[12pt]{article}

\begin{document}
\begin{center}
{\bf GAUSS DECOMPOSITION FOR QUANTUM GROUPS AND DUALITY } \\
\end{center}

{\bf E.V.Damaskinsky}\footnote{
Supported by Russian Foundation for  Fundamental Research, Grant N
95-01-00569-a}{$^,$}\footnote{
Defense Constructing Engineering Institute,  Zacharievskaya st 22, 191194,
St.Petersburg, Russia} {\bf P.P.Kulish}\footnote{
St.Petersburg Branch of Mathematical Institute of Russian Academy of
Science, Fontanka 27, 191011, \hfil \break
St.Petersburg, Russia. e-mail:kulish@lomi.spb.su}
{\bf V.D. Lyakhovsky{$^{1,}$}}\footnote{
St.Petersburg University, Institute of Physics, Petrodvorets, 198904,
Russia. e-mail: lyakhovsky@phim.niif.spb.su} and
{\bf M.A.Sokolov{$^{1,}$}}\footnote{
St.Petersburg Institute of Mashine Building, Poliustrovskii pr 14, 195108,
St Petersburg, Russia. \hfil \break
e-mail: sokol@pmash.spb.su}

\begin{abstract}
The Gauss decomposition of quantum groups and supergroups are considered.
The main attention is paid to the R-matrix formulation of the Gauss
decomposition and its properties as well as its relation to the contraction
procedure. Duality aspects of the Gauss decomposition are also touched. For
clarity of exposition a few simple examples are considered in some details.

Based on the talk, given by the first named author on International
Conference on Symmetry Methods in Physics ({\bf ICSMP-95}) held in JINR
Dubna, Russia, July 10-16, 1995.
\end{abstract}

\section{Introduction}

In this work using elementary algebraic methods and the $R$-matrix approach
\cite{4} we consider the Gauss decomposition of a $q$-matrix $T$
(that is a matrix whose elements are the generators of the considered
quantum group) $T=T_LT_DT_U$
in strictly lower- and upper-triangular matrices (with units on their
diagonals) and a diagonal matrix $T_D$. Such a decomposition of a given
matrix (with non commuting entries) into the product of matrices of the
special type (similar to the Gauss decomposition) is the particular case of
the general factorization problem \cite{12}, which can be considered as a
cornerstone for many constructions of the classical as well as quantum
inverse scattering methods. It should be pointed out that in different
contexts such decompositions can be found (sometimes in non explicit form)
in many papers on quantum deformations (see \cite{80} for more details and
references).

In quite a general framework of the quantum double construction the Gauss
decomposition was considered in \cite{13}, where the universal triangular
objects ${\cal M}^{\pm}$ were defined. Their matrix representations on one
of the factors $(\rho \otimes 1){\cal M}^{\pm}= M^{\pm} $ are the solutions
of the FRT-relations
$$
RM_1^{\pm }M_2^{\pm }=M_2^{\pm }M_1^{\pm }R,\quad
M_1^{+}M_2^{-}=M_2^{-}M_1^{+}.
$$
Their product $M^-M^+$ after the unification of diagonal elements $
(M_{ii}^{-}M_{ii}^{+}=A_{ii})$ gives the $q-$matrix $T$. Usually the new
generators, defined by the Gauss decomposition, have simpler commutation
rules (multiplication) but more complicated expressions for the coproduct.

The simpler structure of the commutation rules for the Gauss decomposition
generators of the quantum groups simplifies \cite{21} the problem of their $
q $-bosonization (that is their realization by the creation and annihilation
operators of the quantum deformed oscillator \cite{37}-\cite{39}). For the
dual objects -- quantum algebras (or quantum deformations of universal
enveloping algebras of the classical Lie algebras) this problem was
considered in \cite{40,41} and for $q$-superalgebras in \cite{42,43}.

The paper is set up as follows. The general properties of the Gauss
decomposition are treated in Sec.2 in the form most suitable for the best
known case of the quantum group $GL_q(n)$, but general methods (after
obvious modifications) are valid also for other quantum groups and even
supergroups. It also concerns the remarks given in this section about, for
example, the contraction procedure and duality considerations valid in other
cases. Some peculiarities of the orthogonal and symplectic quantum groups
and supergroups are discussed in Sec.3. The developed methods are briefly
illustrated in Sec.4 where results of the Gauss decomposition for the
quantum groups $SO_q(3)$, $Sp_q(2)$ and for some supergroups are given. We
refer to the works \cite{80},\cite{81} for more details.

\section{Gauss factorization for the quantum group $GL_q(n)$}

According to the FRT-approach the matrix relation
\begin{equation}
\label{1}RT_1T_2=T_2T_1R,
\end{equation}
encodes quadratic (commutation) relations between quantum group (QG)
generators. The $R$-matrix for $GL_q(n)$ ($A_l$ series) is the
lower-triangular numerical $n^2\times n^2$ -matrix \cite{4}:
\begin{equation}
\label{2.1}R=\sum_{i\neq j}^nE_{ii}\otimes E_{jj}+q\sum_{i=1}^nE_{ii}\otimes
E_{ii}+\lambda \sum_{i<j}^nE_{ji}\otimes E_{ij},
\end{equation}
where $\lambda =q-q^{-1},q\neq 1$ and $E_{kl}$ is the matrix unit, that is
the $n\times n$-matrix with the only nonzero element $(E_{kl})_{kl}=1$. With
the $R$-matrix thus defined the FRT-equation gives the following commutation
rules:
\begin{equation}
\label{2.4}t_{i\,j}t_{ik}=qt_{ik}t_{i\,j},\quad
\,\;t_{ik}t_{l\,j}=t_{l\,j}t_{ik},\quad t_{ik}t_{lk}=qt_{lk}t_{ik},\;\quad
[t_{i\,j},t_{lk}]=\lambda t_{ik}t_{lj},
\end{equation}
where $1\leq j<k\leq n,\;1\leq i<l\leq n$. From these relations it follows
that the algebra $F_q$ has a rich subalgebra structure. In particular, any
four matrix elements $t_{ij},\,t_{kj},\,t_{il},\,t_{kl}\quad (1\leq
i,j,k,l\leq n)$ (in other words, any four elements standing at intersections
of two arbitrary rows and two arbitrary columns) generate a $GL_q(2)$
-subalgebra.

The quantum determinant ($q$-det) is the expression \cite{4}
\begin{equation}
\label{2.5}D_q(T)=\det \nolimits_qT=\sum_\sigma (-q)^{l(\sigma )}t_{1\sigma
(1)}t_{2\sigma (2)}\cdot \ldots \cdot t_{n\sigma (n)},
\end{equation}
where $l(\sigma )$ is the length of a substitution $\sigma \in {\bf S}_n$.
Note that $\det \nolimits_qT$ is the central element for $GL_q(n)$ \cite{4}.
In the $GL_q(2)-$case $q$-det is equal to
$$
{\rm det}_qT=t_{11}t_{22}-qt_{12}t_{21}.
$$
Invertibility of ${\rm det_q}T$ \thinspace is the necessary condition for
endowing the algebra $F_q={\rm Fun}(GL_q(n))$ with a Hopf algebra structure
\cite{4}. We shall suppose this condition to be fulfilled in what follows.

The Gauss decomposition for a $q$-matrix $T$ is
\begin{equation}
\label{2.15}T=T_LT_DT_U=T_LT^{(+)}=T^{(-)}T_U,
\end{equation}
where $T_L=(l_{ik})$ is a strictly lower-triangular matrix (with the units
at the main diagonal: $l_{kk}=1$); $T_D={\rm diag}(A_{kk})$; $T_U=(u_{ik})$
is a strictly upper-triangular matrix ($u_{kk}=1$); $T^{(+)}:=T_DT_U$ and $
T^{(-)}=T_LT_D$.

The Gauss algorithm of triangularization of classical number matrices can be
applied, as well, for quantum groups. Introduce a strictly lower-triangular $
n\times n$-matrix $W_L$, which maps the $q$-matrix $T$ into an
upper-triangular matrix $T^{(+)}$:
\begin{equation}
\label{2.6}W_LT=T^{(+)}.
\end{equation}
As we wish the upper triangularity of the matrix $T^{(+)},$ the condition (
\ref{2.6}) gives the system
\begin{equation}
\label{2.7}w_kT_{(k-1)}=-t_k
\end{equation}
of linear equations for elements of a row $
w_k=(w_{k,1},w_{k,2},...,w_{k,k-1})$ with generators of the algebra $F_q$ as
coefficients. In this system $T_{(k)}$ denotes a matrix obtained from $T$ by
omitting the last $(n-k)$ rows and $(n-k)$ columns. If the $q$-matrix $
T_{(k)}$ is invertible, then the solution of this system has the form
$$
w_k=-t_kT_{(k-1)}^{\,\,-1}.
$$
The elements of the inverse matrix are given by the formulae \cite{4}
\begin{equation}
\label{2.10}(T^{\,-1})_{ij}=(-q)^{i-j}(\det \nolimits_qT)^{-1}M_q(i,j).
\end{equation}
where $M_q(i,j)$ is the $q$-minor, that is the $q$-det of the matrix which
is obtained from $T$ by omitting of the $i$th row and $j$th column. Thus,
all of the $W_L$-matrix elements are uniquely defined. The condition (\ref
{2.6}) allows us to find all the nonzero elements of the matrix $T^{(+)}$.
In particular the diagonal elements (except $(T^{(+)})_{11}=t_{11}$) have
the form
\begin{equation}
\label{2.11}
\begin{array}{c}
(T^{(+)})_{ii}=(\det \nolimits_qT_{(i-1)})^{-1}(\det \nolimits_qT_{(i)}).
\end{array}
\end{equation}
Note that the last relation is a direct $q$-analog of the classical formula.
It is easy to see that the diagonal elements $(T^{(+)})_{ii}$ of the matrix $
T^{(+)}$ are mutually commuting. With the help of the commutation relations
between the elements $t_{ij}$ of the matrix $T$ and their main minors we can
find the commutation relations between the elements $(T^{(+)})_{ij}$ of the
matrix $T^{(+)}$.

The matrix $W_L$ is strictly-lower triangular, so it can be inverted.
Elements of the inverse matrix $(T_L)=(W_L)^{-1}$ are polynomials on the
elements of the matrix $W_L$. The relation (\ref{2.6}) provides the desired
decomposition
\begin{equation}
\label{2.13}T=T_LT^{(+)}.
\end{equation}

In close analogy with the previous case, one can (using an operator $W_U$)
define the decomposition
\begin{equation}
\label{2.13a}T=T^{(-)}T_U.
\end{equation}
Multiplying from the right the relation (\ref{2.13}) by the matrix $W_U$, $
T^{(-)}=TW_U$, we obtain
$$
T^{(-)}=TW_U=T_LT^{(+)}W_U=T_L(T^{(+)}T_U^{\,-1})=T_LT_D,
$$
where $T_D=T^{(+)}T_U^{\,-1}$ is the diagonal matrix with elements of the
form given by the formula (\ref{2.11}). In the same manner we obtain $
T^{(+)}=T_DT_U$. As a result we have the Gauss factorization (\ref{2.15}) of
the $q$-matrix $T=T_LT_DT_U$.

In principle the factorization procedure considered above allows us to get
the expressions for elements of all the matrices participating in the
decomposition in terms of the original generators of  $GL_q(n)$. Moreover,
the procedure defines the explicit form of $\det \nolimits_qT$ and minors
and commutation rules for them. The latter is especially important for the
QGs of series other than $A_n$. Unfortunately, such a factorization
procedure is rather cumbersome. This difficulty can be avoided with the help
of the contraction procedure considered, for example, in \cite{44,PK} for
the case of QGs.

In the fundamental representation the Cartan elements $h_i$ of $gl(n)$ are
realized by the $n\times n$ matrices $h_i=\frac 1{2n}E_{ii}$. They are not
changed by standard quantum deformations. Moreover, they remain primitive,
that is, their coproducts have the form $\Delta (h_i)\equiv H_i=h_i\otimes
1+1\otimes h_i.$ Due to the cocommutativity of Cartan elements ($\Delta
(h_i)=\Delta ^{\prime }(h_i)$ with $\Delta ^{\prime }=P\circ \Delta$) we
have
\begin{equation}
\label{2.16}[R,H_i]=0.
\end{equation}

Let us subject the FRT-relation (\ref{1}) to the similarity transformation
with the matrix
\begin{equation}
\label{2.16a}K_\gamma =\exp (\sum_{i=1}^n\gamma _ih_i)\otimes \exp
(\sum_{i=1}^n\gamma _ih_i)=\exp (\sum_{i=1}^n\gamma _i\Delta (h_i))=\exp
(\sum_{i=1}^n\gamma _iH_i),
\end{equation}
where the numerical coefficients $\gamma _i$ are strictly ordered $\gamma
_1>\gamma _2>\ldots >\gamma _n>0$. In view of the relation (\ref{2.16}),
such a transformation affects only the matrices $T_l,(l=1,2)$:
\begin{equation}
\label{2.16b}T\longrightarrow K_\gamma TK_\gamma ^{\,-1}={\rm diag}
(e^{\gamma _i})T{\rm diag}(e^{-\gamma _i}),\qquad t_{ij}\longmapsto
t_{ij}e^{\gamma _i-\gamma _j}.
\end{equation}
Introduce the two sets of new generators
$$
t_{\;ij}^{(+)}=\left\{
\begin{array}{cc}
t_{ij}e^{\gamma _i-\gamma _j} & i\leq j \\
t_{ij} & i>j
\end{array}
\right. ;\quad t_{\;ij}^{(-)}=\left\{
\begin{array}{cc}
t_{ij} & i<j \\
t_{ij}e^{\gamma _i-\gamma _j} & i\geq j
\end{array}
\right. .
$$
Let $\gamma _i-\gamma _j=\gamma _{ij} \varepsilon$, where $\gamma _{ij}>0$,
if $i<j$ and $\gamma _{ij}<0$, if $i>j$. When $\varepsilon\rightarrow \infty
\;(\varepsilon\rightarrow -\infty )$ in the set $\{t_{\;ij}^{(+)}\}\quad
(\{t_{\;ij}^{(-)}\})$ all the matrix elements with $i>j\,\,\,(i<j)$ vanish.
Thus, we have constructed the homomorphisms of the algebra $F_q$ into the
algebras $F_q^{(\pm )}$, generated by the elements of the upper- and
lower-triangular $q$-matrices. In this case the commutation rules for the
new generators are completely defined by the initial $R$-matrix:
\begin{equation}
\label{2.19}RT_1^{(\pm )}T_2^{(\pm )}=T_2^{(\pm )}T_1^{(\pm )}R.
\end{equation}
Note that for contracted algebras (and their generators) we use the same
notations as in the case of algebras obtained by triangularization
procedure. It can be verified that commutation relations for the Gauss
generators are uniquely defined by the $R$-matrix.

Similar contraction procedure allows us to find the homomorphisms $T^{(\pm
)}\rightarrow A$ of the QGs described by $T^{(\pm )}$ into the group whose
generators are the elements of the diagonal matrix $T_D$. The multiplication
rules for the latter group are also determined by the FRT-relation $
RA_1A_2=A_2A_1R,\qquad (A_{ij}=\delta _{ij}A_{ij}),$ which, in view of the
structure of the $R$-matrix for the quantum group $GL_q(n),$ is equivalent
to the relation
\begin{equation}
\label{2.21}A_1A_2=A_2A_1.
\end{equation}
The relation (\ref{2.21}) means the commutativity of the elements of the
diagonal quantum group.

In the limit $\epsilon \rightarrow \infty $ ($\epsilon \rightarrow -\infty $
) the transformation (\ref{2.16a}) can be treated not only as a
factorization of $F_q$ but also as a contraction procedure for the Hopf
algebras: $F_q\rightarrow F_q^{{\rm contr}}$. In terms of generators $
\left\{ t_{ij}^{(+)}\right\} $ (respectively $\left\{ t_{ij}^{(-)}\right\} $
) this contraction is described by the following transformation:
\begin{equation}
\label{eq20}
\begin{array}{cc}
t_{ij}^{(+)}\rightarrow \left\{
\begin{array}{cc}
t_{ij}^{(+)}; & i\leq j \\
e^{\gamma _{ij}}t_{ij}^{(+)}; & i>j
\end{array}
\right. & {\rm or}\quad t_{ij}^{(-)}\rightarrow \left\{
\begin{array}{cc}
e^{\gamma _{ij}}t_{ij}^{(-)}; & i<j \\
t_{ij}^{(-)}; & i\geq j
\end{array}
\right.
\end{array}
\end{equation}

We have seen that for the set of generators $\left\{ t_{ij}^{(+)}\right\} $
on the factorspace the multiplication is still defined by the $RTT$-equation
with the same $R$ . The limiting transition (\ref{eq20}) allows one to
examine the whole set of generators and it is easily seen that \underline{all
} the structure constants of $F_q$ described by the relations (\ref{2.4})
have the finite limits when $\epsilon $ tends to $+\infty $ (respectively $
-\infty $ ). Obviously the co-structure constants of $F_q$ also have the
finite limit values. Consider, for example, the co-algebra $F_q^{{\rm contr}}
$ for the first type of transformations (\ref{eq20}):
\begin{equation}
\label{eq21}\Delta t_{ij}^{(+)}=\left\{
\begin{array}{ccc}
\sum\limits_{i\leq s\leq j}t_{is}^{(+)}\otimes t_{sj}^{(+)} & {\rm for} &
i<j \\
\sum\limits_{j\leq s\leq i}t_{is}^{(+)}\otimes t_{sj}^{(+)} & {\rm for} & j<i
\end{array}
\right.
\end{equation}
The coproduct here does not mix the upper and lower parts of the $T$-matrix.
Note that while in (\ref{2.19}) the limiting procedure nullifies one of the
triangular parts of the initial $T$-matrix, the relations (\ref{2.4})
rewritten for $t_{ij}^{(+)}$ still describe the multiplication rules for all
the generators of $F_q^{{\rm contr}}$. The contraction does not touch the
compositions of the first three types in (\ref{2.4}). To fix the new
relations of the fourth type let us introduce the grading function $\sigma .$
$$
\sigma \left( t_{ij}^{(+)}\right) =\left\{
\begin{array}{ccc}
+1 & {\rm for} & i<j \\
0 & {\rm for} & i=j \\
-1 & {\rm for} & i>j
\end{array}
\right. .
$$
For $1\leq j<k\leq n$ and $1\leq i<l\leq n$ the last commutator in (\ref{2.4}
) has the following values:
\begin{equation}
\label{eq22}\left[ t_{ij}^{(+)},t_{lk}^{(+)}\right] \!=\!\left\{ \!
\begin{array}{cc}
0 & \left\{
\begin{array}{cc}
{\rm for} & \sigma \left( t_{ij}^{(+)}\right) \neq \sigma \left(
t_{lk}^{(+)}\right)  \\
{\rm for} & \sigma \left( t_{ij}^{(+)}\right) =\sigma \left(
t_{lk}^{(+)}\right) ,\sigma \left( t_{ik}^{(+)}\right) \neq \sigma \left(
t_{lj}^{(+)}\right)
\end{array}
\right.  \\
\lambda t_{ik}^{(+)}t_{lj}^{(+)} & {\rm for}\;\sigma \left(
t_{ij}^{(+)}\right) =\sigma \left( t_{lk}^{(+)}\right) ,\sigma \left(
t_{ik}^{(+)}\right) =\sigma \left( t_{lj}^{(+)}\right)
\end{array}
\right.
\end{equation}
The Hopf algebra $F_q^{{\rm contr}}$ describes a QG. The corresponding
quantum algebra can be defined on the space dual to that of $F_q^{{\rm {contr
}}}$. To obtain the appropriate form of relations for the dual Hopf algebra
let us use the basis $\left\{ l_{ik},u_{ik},A_{kk}\right\} $ for $F_q$ and $
F_q^{{\rm contr}}$. Comparing the left part of (\ref{2.16a}) with the
decomposition (\ref{2.15}) one can see that the transformations (\ref{2.16a}
) and (\ref{2.16b}) have in this basis the same form:
\begin{equation}
\label{eq23}\left\{
\begin{array}{c}
l_{ik}\rightarrow l_{ik}e^{\gamma _{ik}}\equiv l_{ik}^{\left( +\right) } \\
u_{ik}\rightarrow u_{ik}e^{\gamma _{ik}}\equiv u_{ik}^{\left( +\right)
}e^{\gamma _{ik}} \\
A_{kk}\rightarrow A_{kk}\equiv A_{kk}^{\left( +\right) }
\end{array}
\right.
\end{equation}
Let us introduce on the space dual to that of $F_q$ the generators $\left\{
\mu _{i,i+1},\nu _{j+1,j}\alpha _{kk}\right\} ,$
\begin{equation}
\label{eq24}
\begin{array}{ccc}
\left\langle \alpha _{kk},A_{ss}^{\left( +\right) }\right\rangle =\delta
_{ks}; & \left\langle \alpha _{kk},u_{i,i+1}^{\left( +\right) }\right\rangle
=0; & \left\langle \alpha _{kk},l_{j+1,j}^{\left( +\right) }\right\rangle
=0; \\
\left\langle \mu _{i,i+1},A_{ss}^{\left( +\right) }\right\rangle =0; &
\left\langle \mu _{k,k+1},u_{s,s+1}^{\left( +\right) }\right\rangle =\delta
_{ks}; & \left\langle \mu _{i,i+1},l_{j+1,j}^{\left( +\right) }\right\rangle
=0; \\
\left\langle \nu _{i,i+1},A_{ss}^{\left( +\right) }\right\rangle =0; &
\left\langle \nu _{i,i+1},u_{i,i+1}^{\left( +\right) }\right\rangle =0; &
\left\langle \nu _{i,i+1},l_{j+1,j}^{\left( +\right) }\right\rangle =\delta
_{ij}.
\end{array}
\end{equation}
These elements form the Chevalley basis for $U_q(gl(n))$ with the simple
roots $\left\{ \lambda _i\right\} $ and the defining relations \cite{83}

\begin{equation}
\label{eq25}
\begin{array}{c}
\left[ \alpha _{ii},\alpha _{kk}\right] =0; \quad \left[ \alpha _{ii},\mu
_{j,j+1}\right] =\left( \lambda _i,\lambda _j\right) \mu _{j,j+1};\quad
\left[ \alpha _{ii},\nu _{j+1,j}\right] =-\left( \lambda _i,\lambda
_j\right) \nu _{j+1,j};
\end{array}
\end{equation}
\begin{equation}
\label{eq26}\left[ \mu _{j,j+1},\nu _{k+1,k}\right] _{e^{h/2\left\langle
\lambda _j,\lambda _k\right\rangle }} = \delta _{jk}\frac{e^{-h\alpha
_{jj}}-1}{e^{-h}-1};
\end{equation}
\begin{equation}
\label{eq29}\left( {\rm ad}\mu _{j,j+1}\right) ^{1-\left\langle \lambda
_j,\lambda _k\right\rangle }\mu _{k,k+1}=0; \qquad \left( {\rm ad}\nu
_{j,j+1}\right) ^{1-\left\langle \lambda _j,\lambda _k\right\rangle }\nu
_{k,k+1}=0;
\end{equation}
\begin{equation}
\label{eq30}
\begin{array}{c}
\Delta \alpha _{ii}=\alpha _{ii}\otimes 1+1\otimes \alpha _{ii}; \\
\Delta \mu _{j,j+1}=\mu _{j,j+1}\otimes 1+e^{-\frac h2\alpha_{jj}}\otimes
\mu _{j,j+1};\qquad \Delta \nu _{j,j+1}=\nu _{j,j+1}\otimes 1+e^{-\frac
h2\alpha _{jj}}\otimes \nu _{j,j+1};
\end{array}
\end{equation}
The transformation
\begin{equation}
\label{eq33}\nu _{i,j}\rightarrow e^{-\gamma _{ij}}\nu _{i,j}, \qquad \mu
_{i,j}\rightarrow \mu _{i,j}, \qquad \alpha _{ii}\rightarrow \alpha _{ii}
\end{equation}
is dual to (\ref{eq23}) and describes a contraction $U_q(gl(n))\rightarrow
U_q^{{\rm contr}}$ when $\epsilon \rightarrow \infty $ . In this limiting
procedure the duality is preserved. The obtained Hopf algebra $U_q^{{\rm
contr}}$ is dual to $F_q^{{\rm contr}}$. The contraction changes only one of
the defining relations: instead of (\ref{eq26}) one obtains
\begin{equation}
\label{eq34}\left[ \mu _{j,j+1},\nu _{k+1,k}\right] _{e^{h/2\left\langle
\lambda _j,\lambda _k\right\rangle }}=0;
\end{equation}
Thus $U_q^{{\rm contr}}$ is the quantization of a Lie algebra $gl(n)^{{\rm
contr}}$. The latter has the same Borel subalgebras $b^{+}$and $b^{-}$ as
the original $gl(n)$ but here the corresponding $n^{+}$ and $n^{-}$
subalgebras commute:
\begin{equation}
\label{eq35}\left[ n^{+},n^{-}\right] =0
\end{equation}
This commutativity is in total accordance with the ''separation'' of
coproducts in (\ref{eq21}).

We have thus demonstrated that the contraction described by the relations (
\ref{2.16a}),(\ref{2.16b}) leads to the QG $F_q^{{\rm contr}}$ which is the
quantized algebra of functions over the Lie group $GL^{{\rm contr}}(n,$C$)$
with the Lie algebra $gl^{{\rm contr}}(n,$C$)$ defined by the relations (\ref
{eq25}), (\ref{eq29}-\ref{eq30}) and (\ref{eq34}).

The contraction of the type (\ref{eq33}) can be performed for an arbitrary
simple Lie algebra. It is sufficient to multiply every element $x_\lambda $
of $b^{+}$ (or $b^{-}$) by $\epsilon ^{-\sum m_i}$ (respectively $\epsilon
^{\sum m_i}$) and go to the limit $\epsilon \rightarrow \infty $. Here $m_i$
are the coordinates of the root $\lambda $ in terms of simple roots $\left\{
\lambda _i\right\} $. Such a contraction always exists and forces
subalgebras $n^{+}$and $n^{-}$ to commute.

Considering contractions performed by the operators $K_{1\gamma }=\exp
(\sum_{i=1}^n\gamma _ih_i)\otimes 1$ and $K_{2\gamma }=1\otimes \exp
(\sum_{i=1}^n\gamma _ih_i),$ we obtain the following relations (where $R_D$
is the diagonal part of the $R$-matrix).
\begin{equation}
\label{2.22}R_DA_1T_{\;2}^{(-)}=T_{\;2}^{(-)}A_1R_D;\quad
R_DT_{\;1}^{(+)}A_2=A_2T_{\;1}^{(+)}R_D;\quad
R_DT_{\;1}^{(+)}T_{\;2}^{(-)}=T_{\;2}^{(-)}T_{\;1}^{(+)}R_D.
\end{equation}
Other commutation rules for the generators of the QGs $T_U$ and $T_L$ can be
obtained from the relations (\ref{2.15},\ref{2.19}-\ref{2.22}):
\begin{equation}
\label{2.25}
\begin{array}{cc}
RR_DT_{U1}R_D^{\;-1}T_{U2}=R_DT_{U2}R_D^{\;-1}T_{U1}R,\, &
RT_{L1}R_DT_{L2}R_D^{\;-1}=T_{L2}R_DT_{L1}R_D^{\;-1}R, \\
R_DT_{D1}T_{L2}=T_{L2}T_{D1}R_D, & R_DT_{U1}T_{D2}=T_{D2}T_{U1}R_D,
\end{array}
\end{equation}
(and similar equalities with interchange $1\leftrightarrow 2$). If we
separate the $T_U$ and $T_L$ parts from $T^{(\pm )}$ in (\ref{2.22}) and
take into account relations (\ref{2.25}) we get
\begin{equation}
\label{2.27a}T_{U1}T_{L2}=T_{L2}T_{U1}.
\end{equation}
Thus, all the elements of the $q$-matrix $T_U$ commute with every element of
the $q$-matrix $T_L$. Because of the relation $[R,R_D]=0$, which is valid
for $GL_q(n)$, $Sp_q(n)$ and $SO_q(2n)$ (but not for $SO_q(2n+1)$!), we can
rewrite the first two equations in (\ref{2.25}) in the form of the
reflection equation
$$
RT_{U1}R_D^{\;-1}T_{U2}=T_{U2}R_D^{\;-1}T_{U1}R,\quad
RT_{L1}R_DT_{L2}=T_{L2}R_DT_{L1}R.
$$

The relations (\ref{2.15},\ref{2.19},\ref{2.22},\ref{2.25},\ref{2.27a})
together with easily deduced ones
\begin{equation}
\label{2.26}R_DT_1^{(+)}T_{L2}=T_{L2}R_DT_1^{(+)};\qquad T_2^{(-)}{R_D}
^{-1}T_{L1}=T_{L1}T_2^{(-)}{R_D}^{-1},
\end{equation}
supply us with the full sets of commutation rules imposed on the elements of
the matrices $\{T_L,T^{(+)}\},$ $\{T^{(-)},T_U\}$ and $\{T_L,T_D,T_U\}$.
This allows to consider each of these sets of elements as the new basis of
generators for $GL_q(n)$. The basis \{$T_L,T_D,T_U$\}, which has the maximal
number (for $N\geq 3$) of mutually commuting elements is a particularly
convenient. Diagonal elements of $T_D$ play the role of Cartan generators of
a Lie algebra. $Q$-det (\ref{2.5}) is the central element of $GL_q(n)$ and,
as in the case of numerical matrices, can be written as the product of the
diagonal matrix elements
\begin{equation}
\label{2.28-1}\det \nolimits_qT=\prod_{i=1}^n(T_D)_{ii}.
\end{equation}
Of course, $\det \nolimits_qT$ commutes with all the new generators just as
well as with the original ones.

In conclusion of the Sec. we note that further decomposition of $T_L$ and $
T_U$ is possible. Such a procedure leads to a new basis in $GL_q(n)$ with
Weyl-type commutation relations for its elements \cite{84}.

\section{Orthogonal and symplectic quantum groups, quantum supergoups}

Commutation relations for generators of the orthogonal and symplectic QGs
are determined not only by FRT-equation with the relevant $R$-matrix, but
also by supplementary conditions \cite{4}
\begin{equation}
\label{3.2}TCT^tC^{-1}=1,\qquad CT^tC^{-1}T=1,
\end{equation}
which are quantum analog of the known conditions for matrices of the
orthogonal and symplectic Lie groups in the defining representation. In (\ref
{3.2}) $T^t$ is the transposition of $T$, $C$ is a numerical matrix \cite{4}
. Among the supplementary conditions encoded by (\ref{3.2}) only those that
have the nonzero right hand side are linearly independent with respect to
the commutation relations defined by the FRT-equation (\ref{1}).

Similarly to the case of $GL_q(n)$ in the series $B_n,C_n$ and $D_n$ (with
the only exception of $B_1$) the corresponding quantized algebras of
coordinate functions $F_q$ have rich subalgebra structures. For example, a $
GL_q(k)$-subalgebra is generated by the elements of the $q$-matrix $T\,$
located at the intersections of rows with indices $(i_1,i_2,..,i_k)$ and
columns with indices $(j_1,j_2,..,j_k)$ (when $k\leq n$ and no pairs of the
type $(i_m,i_l)$ with $i_m=(i_l)^{\prime }=N+1-{i_l}$ occur in the set $
\{i_p\}$ and also in $\{j_p\}$). The products of the generators belonging to
the different $GL_q(k)$-subalgebras have much more complicated form.

The Gauss factorization of orthogonal, symplectic groups and supergroups can
be carried out by the method described in section 2 for $GL_q(n).$ The
corresponding commutation relations are governed by the same formulae ((\ref
{2.15},\ref{2.19},\ref{2.22},\ref{2.25},\ref{2.27a},(\ref{2.26}). One must
take into consideration that the supplementary conditions yield some new
relations. For instance, $T_D$-matrix generators satisfy the relations
\begin{equation}
\label{3.14}(T_D)_{ii}(T_D)_{i^{\prime }i^{\prime }}=1,\quad 1\leq i\leq N.
\end{equation}

Using the multiplication rules one can prove that the supplementary
conditions for the matrices $T^{(\pm )}$ have the same form (\ref{3.2}) as
for $T$. It should be noted, that with respect to the elements of $T_L$ (or $
T_U$) these conditions contain not only quadratic but also linear terms. It
allows us to exclude dependent generators from the generator list. The
number of remaining generators is equal to the dimension of the
corresponding Lie group. We shall illustrate such a reduction in the next
section.

\section{Examples}

In this Sec. some details of the Gauss factorization of the symplectic and
orthogonal QGs are considered using SO$_q$(3) and $Sp_q(2)$ as examples. The
last part of this Sec. is devoted to the exposition of some simplest
examples of quantum supergroups. The $R$-matrices for these groups have more
nonzero matrix elements in comparison with the $R$-matrix for $GL_q(n)$ of
the same rank. The complicated structure of the $R$-matrix causes additional
complexification of the commutation relations (even for the already
mentioned simplest cases the corresponding list of relations is rather long
and cannot be presented here). Gauss decomposition changes this situation
drastically. It provides the reduced list of generators, containing only
independent ones, with rather simple commutation rules.

\subsection{Quantum group $B_1\sim SO_q(3)$}

Let us apply the Gauss factorization procedure to the $q$-matrix of  $SO_q(3)
$
\begin{equation}
\label{4.1}T=T_LT_DT_U=\left(
\begin{array}{ccc}
1 & 0 & 0 \\
l_{21} & 1 & 0 \\
l_{31} & l_{32} & 1
\end{array}
\right) \left(
\begin{array}{ccc}
A_{11} & 0 & 0 \\
0 & A_{22} & 0 \\
0 &  & A_{33}
\end{array}
\right) \left(
\begin{array}{ccc}
1 & u_{12} & u_{13} \\
0 & 1 & u_{23} \\
0 & 0 & 1
\end{array}
\right)
\end{equation}
After cumbersome computations one finds that here all the matrix elements
can be described in terms of three independent generators $u,\,A$ and $z$
(just as in the classical Lie case):
\begin{equation}
\label{4.2}
\begin{array}{c}
A_{11}=A,\qquad A_{22}=1,\qquad A_{33}=A^{-1}; \\
\begin{array}{lll}
u_{12}=t_{11}^{\;-1}t_{12}=u,\; & u_{13}=t_{11}^{\;-1}t_{13}=-\mu
^{-1}u^2,\; & u_{23}=-q^{1/2}u; \\
l_{21}=t_{21}t_{11}^{\;-1}=z,\; & l_{31}=t_{31}t_{11}^{\;-1}=-\mu
^{-1}z^2,\; & l_{32}=-q^{-1/2}z.
\end{array}
\end{array}
\end{equation}
This gives
\begin{equation}
\label{4.3}
\begin{array}{c}
T=T_LT_DT_U=\left(
\begin{array}{ccc}
A & Au & -\mu ^{-1}Au^2 \\
zA & I+zAu & -\mu ^{-1}zAu^2-q^{1/2}u \\
-\mu ^{-1}z^2A & -\mu ^{-1}z^2Au-q^{-1/2}z & \mu ^{-2}z^2Au^2+zu+A^{-1}
\end{array}
\right)
\end{array}
\end{equation}

It is easy to see, that these new generators are subject to the following
Weyl-like commutation rules
\begin{equation}
\label{4.4}Au=quA,\quad Az=qzA,\quad uz=zu.
\end{equation}

As a simple application of these results let us show that the $q$-det for
the QG is equal to unity. We propose the following notations for minors of
the quantum matrix $T$
\begin{equation}
\label{4.5}\Delta \left[
\textstyle{
{{i\;j}\atop {k\;l}};q^\alpha }\right]
\equiv t_{ik}t_{jl}-q^\alpha t_{il}t_{jk},
\end{equation}
Then the $q$-det can be written in a simple form
\begin{equation}
\label{4.6}
\begin{array}{c}
\det \nolimits_qT=t_{11}\Delta \left[ \textstyle{{{2\;3}\atop {2\;3}};q}
\right] -t_{12}\Delta \left[ \textstyle{{{2\;3}\atop {1\;3}};q^2}\right]
+qt_{13}\Delta \left[ \textstyle{{{2\;3}\atop {1\;2}};q}\right]
\end{array}
\end{equation}
In \cite{11} the same expression was justified by geometric considerations.
If we rewrite it in terms of new generators we easily get
\begin{equation}
\label{4.7}\det \nolimits_qT=1.
\end{equation}
This result agrees with the general statement made for QGs of classical
types(see Subsec. 2.3 above): the $q$-det is equal to the product of
diagonal elements of the matrix $T_D.$

The single matrix element can also be written as $t_{ij}= \Delta \left[
\textstyle{{i}\atop {j}}\right] $. In these terms the nontrivial
elements of $T_L$ and $T_U$ are the ratios of $q$-minors in complete analogy
with classical (commutative) situation. For example,
\begin{equation}
\label{4.8}T_U=\left(
\begin{array}{ccc}
1 & \Delta ^{-1}\left[ \textstyle{{1}\atop {1}}\right] \Delta \left[
\textstyle{{1}\atop {2}}\right] & \Delta ^{-1}\left[
\textstyle{{1}\atop {1}}\right] \Delta \left[\textstyle{{1}\atop {3}}
\right] \\ 0 & 1 & \vphantom{\biggl[} \Delta ^{-1}\left[ \textstyle
{{{2\;1}\atop {1\;2}};q}\right] \Delta \left[
\textstyle{{{2\;1}\atop {1\;3}};q^2} \right] \\ 0 & 0 & 1 \end{array}
\right)
\end{equation}

\subsection{Quantum group $C_2\sim Sp_q(2)$}

To realize the Gauss algorithm we must find the operators $W_L$ and $W_U$
(see Sec.2). Solving, for example, the system of equations (\ref{2.6}) for $
W_L$ one has (in notation introduced in (\ref{4.5}))
\begin{equation}
\label{5.7}
\begin{array}{c}
w_{21}^L=-t_{21}\Delta _q^{\,-1}[
\textstyle{{1}\atop {1}}];\quad w_{31}^L=(\Delta _q[\textstyle
{{{2\,3}\atop {1\,3}};q}]- \lambda \Delta
_q[\textstyle{{{1\,4}\atop {1\,2}};q}]) \Delta
_q^{\,-1}[\textstyle{{1\,2}\atop {1\,2}}];\quad w_{32}^L=-\Delta
_q[\textstyle{{{1\,3}\atop {1\,2}};q}] \Delta _q^{\,-1}[ \textstyle
{{{1\,2}\atop {1\,2}};q}]; \\ \hfill \\ w_{41}^L=-q^2t_{41} \Delta
_q^{\,-1}[ \textstyle{{1}\atop {1}}];\quad w_{42}^L= -q^2t_{31}\Delta
_q^{\,-1}[ \textstyle{{1}\atop {1}}];\quad w_{43}^L= t_{21}\Delta
_q^{\,-1}[ \textstyle{{1}\atop {1}}]; \end{array} \end{equation} Thus,
we obtain $$ \begin{array}{c} T^{(+)}=\left( \begin{array}{cccc} \Delta
_q[\textstyle{{1}\atop {1}}] & t_{12} & t_{13} & t_{14} \\ 0 & \Delta
_q^{\,-1}[\textstyle{{1}\atop {1}}] \Delta _q[\textstyle{
{{1\,2}\atop {1\,2}};q}] & \Delta _q^{\,-1}[\textstyle{{1}\atop {1}}]
\Delta _q[\textstyle{{{1\,2}\atop {1\,3}};q}] & \Delta _q^{\,-1}[
\textstyle{{1}\atop {1}}](t_{11}t_{24}-q^2t_{14}t_{21}) \\ 0 & 0 &
\Delta_q^{\,-1}[\textstyle{{{1\,2}\atop {1\,2}};q}] \Delta _q[\textstyle{
{1}\atop {1}}] & -\Delta_q^{\,-1}[
\textstyle{{1}\atop {1}}]t_{12} \\ 0 & 0 & 0 & \Delta _q^{\,-1}[
\textstyle{{1}\atop {1}}]
\end{array}
\right)
\end{array}
$$

Here in close analogy with the commutative case, the diagonal elements are
the ratios \\ $
\begin{array}{c}
\Delta _q^{\,-1}[\textstyle{{1,2,...,k-1}\atop {1,2,...,k-1}}]\Delta _q[
\textstyle{{1,2,...,k}\atop {1,2,...,k}}]
\end{array}
$ of the diagonal minors. In the matrix given above these ratios are further
simplified using particular properties of the symplectic groups. Thus
evaluating the matrix elements of $T^{(+)}$ on the places related to the
diagonal $GL_q$-minors one fins the expressions of the form
\begin{equation}
\label{5.9}
\begin{array}{c}
D_q^{sp}[\textstyle{{1,2,...,k}\atop {1,2,...,k}}]=\sum\limits_\sigma
(-q)^{l(\sigma )}q^{l^{\prime }(\sigma )}t_{1,\sigma (1)}t_{2,\sigma
(2)}\cdot ...\cdot t_{k,\sigma (k)},
\end{array}
\end{equation}
with additional factor $q^{l^{\prime }(\sigma )}$ in comparison with (\ref
{2.5}). Here $l^{\prime }(\sigma )$ is the number of transpositions of
''specific'' elements (the transposition index). For example, $l^{\prime
}(1,2,4,3)=0$, but $l^{\prime }(1,3,2,4)=1$ (2 and $3=2^{\prime }$ are
transposed). For $i=3,4$ the following simplifications become possible in (
\ref{5.9}),
\begin{equation}
\label{5.10}
\begin{array}{c}
D_q^{sp}[\textstyle{{1,2,3}\atop {1,2,3}}]=D_q^{sp}[\textstyle{{1}\atop
{1}}]\qquad D_q^{sp}(T)=D_q^{sp}[\textstyle{{1,2,3,4}\atop {1,2,3,4}}
]=1.
\end{array}
\end{equation}
So, in the symplectic case it is natural to define the $q$-det by the
formula (\ref{5.9}). This definition is in agreement with the one presented
in \cite{11} with the different argumentation.

The non zero elements of the matrix $W_L^{\;-1}=T_L$ have the form (remind
that $T_L=(l_{ij}), T_D=(A_{ii}), T_U=(u_{ij})$)
\begin{equation}
\label{5.11}
\begin{array}{c}
(T_L)_{21}=-w_{21};\quad (T_L)_{32}=-w_{32};\quad (T_L)_{43}=-w_{43}; \\
\hfill \\ (T_L)_{31}=w_{32}w_{21}-w_{31}=t_{31}D_q^{\,-1}[
\textstyle{{1}\atop {1}}]; \\ \hfill \\
(T_L)_{41}=-w_{43}w_{32}w_{21}+w_{43}w_{31}+w_{42}w_{21}-w_{41}=t_{41}
D_q^{\,-1}[
\textstyle{{1}\atop {1}}]; \\ \hfill \\
(T_L)_{21}=w_{43}w_{32}-w_{42}=q^{-1} D_q[ \textstyle{{1,4}\atop {1,2}}
]D_q^{\,-1}[ \textstyle{{1,2}\atop {1,2}}];
\end{array}
\end{equation}
Using definition $T^{(+)}=T_DT_U$ matrix elements of $T^{(+)}$ can be
obtained. Analogous procedure based on the matrix $W_U$ leads, naturally, to
the same results.

The obtained formulas allow us to find commutation relations for the new
basic generators induced by the Gauss decomposition. The final form of these
relations is
\begin{equation}
\label{5.12}
\begin{array}{c}
\left\{
\begin{array}{l}
\lbrack A_{kk},A_{jj}]=0 \\
\lbrack u_{kl},l_{ij}]=0
\end{array}
\right. \qquad \left\{
\begin{array}{rcl}
A_ml_{ij} & = & q^{(\delta _{mj}-\delta _{mj^{\prime }}-\delta _{mi}+\delta
_{mi^{\prime }})}l_{ij}A_m; \\
A_mu_{ij} & = & q^{(\delta _{im}\,-\delta _{im^{\prime }}-\delta
_{jm}+\delta _{jmi^{\prime }})}u_{ij}A_m;
\end{array}
\right. \\
\smallskip \\ {\left\{
\begin{array}{rcl}
l_{21}l_{31} & = & q^2l_{31}l_{21}+q\lambda l_{41}; \\
l_{32}l_{21} & = & q^2l_{21}l_{32}-(q^4-1)l_{31};
\hspace{1cm} \\ l_{31}l_{32} & = & q^2l_{32}l_{31}; \\
\lbrack l_{41},l_{ij}] & = & 0;
\end{array}
\right. }\qquad {\left\{
\begin{array}{rcl}
u_{12}u_{13} & = & q^2u_{13}u_{12}+q\lambda u_{14}; \\
u_{23}u_{12} & = & q^2u_{12}u_{23}-(q^2-q^{-2})u_{13};
\hspace{.6cm} \\ u_{13}u_{23} & = & q^2u_{23}u_{13}; \\
\lbrack u_{14},u_{kl}] & = & 0;
\end{array}
\right. }
\end{array}
\end{equation}

The supplementary conditions cause the following constrains on the set of
new generators
\begin{equation}
\label{5.13}
\begin{array}{c}
A_{33}=A_{22}^{\,\,-1};\qquad \quad \qquad A_{44}=A_{11}^{\,\,-1}; \\
l_{42}=q^2l_{31}-l_{21}l_{32};\quad \quad \quad l_{43}=-l_{21}; \\
u_{24}=q^2(u_{13}-u_{12}u_{23});\quad u_{34}=-u_{12};
\end{array}
\end{equation}
Hence, the number of independent generators decreases to ten, i.e. to the
dimension of the corresponding classical Lie group.

\subsection{Quantum supergroups}

The FRT-relation for quantum supergroups has the same form as in (\ref{1}),
but the matrix tensor product includes additional sign factors (${\pm }$)
related to ${\bf Z}_2$-grading \cite{51}. $Z_2$-graded vector space
(superspace) decomposes into the direct sum of subspaces $V_0\oplus V_1 $ of
even and odd vectors. The parity function ($p(v)=0$ at $v\in V_0$ and $p(w)=1
$ at $w\in V_1$) is defined on them. As a rule, a vector basis with definite
parity $p(v_i)=p(i)=0,1$ is used. In this basis the row and column parities
are introduced in the matrix space {\rm End}$(V)$. The tensor product of two
even matrices $F$, $G$ $(p(F_{ij})=p(i)+p(j))$ has the following signs \cite
{51}
\begin{equation}
\label{6.1}\left( F\otimes G\right)
_{ij;kl}=(-1)^{p(j)(p(i)+p(k))}F_{ik}G_{jl}.
\end{equation}
Due to this prescription $T_2=I\otimes T$ has the same block-diagonal form
as in the usual (non super) case while $T_1=T\otimes I$ includes the
additional sign factor $(-1)$ for odd elements standing at odd rows of
blocks. For the $GL_q(n|m)$ quantum supergroup the $R$-matrix structure is
the same as for the $GL_q(n+m)$ but at odd-odd rows $q$ is changed by $
q^{-1} $
\begin{equation}
\label{6.2}R=\sum_{i,j}\left( 1-\delta _{ij}(1-q^{1-2p(i))})\right)
e_{ii}\otimes e_{jj}+\lambda \sum_{i>j}e_{ij}\otimes e_{ji}.
\end{equation}
Let us remind that the tensor product notation in (\ref{6.2}) refers to the
graded matrices.

The same contraction procedure as in Sec.2 results in homomorphisms of $
T=(t_{ij})$ onto $T^{(\pm )},$ $T_D$ matrices and leads to the corresponding
$R$-matrix relations. Let us give some of them emphasizing the peculiarities
of the supergroup case.

Using the $R$-matrix block structure in the relations
\begin{equation}
\label{6.4}RT_{\quad 1}^{(\pm )}T_{\quad 2}^{(\pm )}=T_{\quad 2}^{(\pm
)}T_{\quad 1}^{(\pm )}R,
\end{equation}
we can find commutation rules for the diagonal elements,
\begin{equation}
\label{6.5}R_D(T_D)_1T_{\quad 2}^{(-)}=T_{\quad 2}^{(-)}(T_D)_1R_D,\quad
R_DT_{\quad 2}^{(+)}(T_D)_1=(T_D)_1T_{\quad 2}^{(+)}R_D.
\end{equation}
As above for the mutually commutative elements $A_{ii}:\,\,T_D={\rm diag}
(A_{11},A_{22},...)$ one has
\begin{equation}
\label{6.6}A_{ii}T^{(\pm )}A_{ii}^{-1}=(R_D)_{\quad ii}^{\pm 1}T^{(\pm
)}(R_D)_{\quad ii}^{\mp 1}.
\end{equation}
However, the diagonal block structure of $R_D$ is different here. As a
consequence the $GL_q(n|m)$ central element is the ratio of the two products
corresponding to even and odd rows
\begin{equation}
\label{6.7}{\rm s-det}_qT=\left( \frac{\prod_{i=1}^nA_{ii}}
{\prod_{k=1}^mA_{n+kn+k}}\right) .
\end{equation}
It is naturally to call this expression the quantum superdeterminant ($q$
-Berezinian).

For the supergroup $GL_q(1|1)$ the commutation relations of the $q$-matrix
elements $T=\left( \textstyle{{a\,\beta }\atop {\gamma \,d}}\right) $
have the form \cite{59}
\begin{equation}
\label{6.15}
\begin{array}{ccl}
a\beta =q\beta a, & \beta d=q^{-1}d\beta , & \beta \gamma =-\gamma \beta ,
\\
a\gamma =q\gamma a, & \gamma d=q^{-1}d\gamma , & ad=da+\lambda \gamma \beta
,
\end{array}
\,\,\,{\beta }^2={\gamma }^2=0.
\end{equation}
Here we used the Greek letters for odd (nilpotent) generators. The Gauss
decomposition gives
\begin{equation}
\label{6.8}T=\left(
\begin{array}{cc}
a & \beta \\
\gamma & d
\end{array}
\right) =\left(
\begin{array}{cc}
1 & 0 \\
\varsigma & 1
\end{array}
\right) \left(
\begin{array}{cc}
A & 0 \\
0 & B
\end{array}
\right) \left(
\begin{array}{cc}
1 & \psi \\
0 & 1
\end{array}
\right) ,
\end{equation}
\begin{equation}
\label{6.9}A=a,\quad \psi =A^{-1}\beta ,\quad \varsigma =\gamma A^{-1},\quad
B=d-\gamma A^{-1}\beta .
\end{equation}
The relations
\begin{equation}
\label{6.10}
\begin{array}{c}
\lbrack A,B]=0,\quad A\psi =q\psi A,\quad A\varsigma =q\varsigma A,\quad
\psi ^2=\varsigma ^2=0, \\
\psi \varsigma +\varsigma \psi =0,\quad B\psi =q\psi B,\quad B\varsigma
=q\varsigma B
\end{array}
\end{equation}
cause the centrality of the superdeterminant in $GL_q(1|1)$ \cite{59}
\begin{equation}
\label{6.11}{\rm s}-\det \nolimits_qT=AB^{-1}=a^2(ad-q\gamma \beta
)^{-1}=a/(d-\gamma a^{-1}\beta ).
\end{equation}

In the $GL_q(2|1)$ case the $q$-matrix of generators has the form
\begin{equation}
\label{6.12}T= \left(
\begin{array}{ccc}
a & b & \alpha \\
c & d & \beta \\
\gamma & \delta & f
\end{array}
\right) = \left(
\begin{array}{cc}
M & \Psi \\
\Phi & f
\end{array}
\right) .
\end{equation}
The even $M$-matrix elements form $GL_q(2)$ subgroup with the commutation
rules (\ref{2.5}). The elements of each ($2\times 2$) submatrix with even
generators at its diagonal form $GL_q(1|1)$ supersubgroup with (\ref{6.15})
- type commutation relations. The remaining multiplications look like
follows
$$
\begin{array}{ccccc}
\alpha\beta=-q^{-1}\beta\alpha,\quad & c\alpha=\alpha c,\quad & b\gamma=
\gamma b, & a\beta=\beta a+\lambda c\alpha,\quad & b\beta=\beta b+ \lambda
d\alpha, \\
\gamma\delta=-q^{-1}\delta\gamma,\quad & d\alpha=\alpha d,\quad &
d\gamma=\gamma d, & a\delta=\delta a+\lambda \gamma b,\quad & c\delta=\delta
c+\lambda \gamma d.
\end{array}
$$

The new generators produced by the Gauss decomposition
\begin{equation}
\label{6.18}T= \left(
\begin{array}{ccc}
a & b & \alpha \\
c & d & \beta \\
\gamma & \delta & f
\end{array}
\right) = \left(
\begin{array}{ccc}
1 & 0 & 0 \\
u & 1 & 0 \\
v & w & 1
\end{array}
\right) \left(
\begin{array}{ccc}
A & 0 & 0 \\
0 & B & 0 \\
0 & 0 & C
\end{array}
\right) \left(
\begin{array}{ccc}
1 & x & y \\
0 & 1 & z \\
0 & 0 & 1
\end{array}
\right)
\end{equation}
have the following commutation rules
$$
\begin{array}{ccl}
Ax=qxA,\quad & Ay=qyA,\quad & Az=zA, \\
Au=quA,\quad & Av=qvA,\quad & Aw=wA,
\end{array}
$$
$$
\begin{array}{ccl}
Bx=q^{-1}xB,\quad & By=yB,\quad & Bz=qzB, \\
Bu=q^{-1}uB,\quad & Bv=vB,\quad & Bw=qwB,
\end{array}
$$
$$
\begin{array}{ccl}
Cx=xC,\quad & Cy=qyC,\quad & Cz=qzC, \\
Cu=uC,\quad & Cv=qvC,\quad & Cw=qwC,
\end{array}
$$
$$
[A,B]=[A,C]=[B,C]=0,\qquad y^2=z^2=v^2=w^2=0
$$
$$
\begin{array}{ccl}
xy=qyx,\quad & yz=-q^{-1}zy,\quad & qxz-zx=\lambda y, \\
uv=qvu,\quad & vw=-q^{-1}wv,\quad & uw-qw^{-1}=\lambda v,
\end{array}
$$
$$
\begin{array}{c}
[x,u]=[x,v]=[x,w]=0,
\quad[u,x]=[u,x]=[u,y]=[u,z]=0, \\ yv+vy=0=yw+wy,\qquad zv+vz=0=zw+wz.
\end{array}
$$

The superdeterminant
$$
{\rm s-}det_qT=ABC^{-1}=det_qM/C
$$
is a central element. The latter expression follows from the block Gauss
decomposition of (\ref{6.12}). In particular for the $GL_q(m|n)$ matrix $T$
in the block form one has (cf. \cite{64})
$$
s-det_qT=det_qA/det_q(D-CA^{-1}B)
$$
which is formally the standard expression.

Generalization of the above results to $GL_q(m|n)$ and other quantum
supergroups looks rather straightforward. Although, as usual, quantum
supergroup $OSp_q(1|2)$ with the rank one has its own peculiarities. In this
case the $q$-matrix $T$ has three independent generators while the
undeformed supergroup has five. The $T$-matrix for this quantum supergroup
in the fundamental co-representation has the dimension three and three
diagonals of the 9x9 $R$-matrix diagonal blocks are $(q,1,1/q)\;(1,1,1) $
and $(1/q,1,q)$. Hence, the diagonal elements $T_D={\rm diag}(A,\,B,\,C)$ of
its Gauss decomposition give rise to the central elements $AC=CA$ and $B.$
The quantum super-determinant is $s-det_qT=AC/B=B$ due to the supplementary
quadratic relation $T^{st}CT=\gamma C$ with $\gamma =B^2=AC=CA$ as in the
case of the orthogonal and symplectic series. The lower-triangular matrix $
T_L$ (as well as upper-triangular $T_U$) has only one independent generator:
$$
T_U=\left(
\begin{array}{ccc}
1 & x & x^2/\omega \\
0 & 1 & -x/q^{1/2} \\
0 & 0 & 1
\end{array}
\right) ,\qquad T_L=\left(
\begin{array}{ccc}
1 & 0 & 0 \\
u & 1 & 0 \\
u^2/\omega & q^{1/2}u & 1
\end{array}
\right) ,
$$
here $\omega =q^{1/2}-q^{-1/2}$. This fact is also reflected in the
structure of the universal $T$-matrix : ${\cal T}
=E_q^{(s)}(V_{-}u)exp(2Ha)E_{1/q}^{(s)}(V_{+}x)$, where $x,\,u$ and $exp(a)=A
$ are generators for the $q$-super-group $OSp_q(1|2)$, while $V_{-}\;V_{+}$
and $H$ are the generators of the dual Hopf super-algebra $osp_q(1|2)$ and $
E_q^{(s)}(t)$ is the $q$-exponent (see also \cite{16}, \cite{84})

\section{Conclusion}

In this work we considered the Gauss decomposition of the quantum groups
related to the classical Lie groups and supergroups by the elementary linear
algebra and $R$-matrix methods. The Gauss factorization yields a new basis
for these groups which is sometimes more convenient than the original one.
Most of the relations for the Gauss generators are written in the $R-$matrix
form. These commutation relations are simpler then the original rules. This
is especially evident in the symplectic and orthogonal cases. In terms of
new generators the supplementary conditions characteristic to the
definitions of $B,C,D$-series of quantum groups allow to extract the
independent ones. Their number is equal to the dimension of the
corresponding classical group. The Gauss factorization leads naturally to
appearance of $q$-analogs of such classical notions as determinants,
superdeterminants and minors. We also want to stress that the new basis is
helpful for studies in quantum group representation theory. In particular as
demonstrated in \cite{21} it simplifies the $q$-bosonization problem. As was
pointed out in the Introduction it looks that almost any relation and/or
statement for standard matrices being appropriately ''$q-$deformed'' is
valid for $q-$matrices.

{\bf Acknowledgement}.

The authors would like to thank the organizers of the conference, especially
G. Pogosyan and P. Winternitz, for warm hospitality rendered to one of them
(E.V.D) at Dubna.

The investigations presented in this article are supported by the Russian
Foundation for Fundamental Research under the Grant N 95-01-00569-a.

\end{document}